\documentclass[english]{article}
\usepackage{amssymb}
\usepackage{amsmath}
\usepackage{babel}

\def\thebibliography#1{\section*{References}\list
  {[\arabic{enumi}]}{\settowidth\labelwidth{[#1]}\leftmargin\labelwidth
    \advance\leftmargin\labelsep
    \usecounter{enumi}}
    \def\newblock{\hskip .11em plus .33em minus -.07em}
    \sloppy
    \sfcode`\.=1000\relax}
\newcommand{\refbook}[3]{{\sc #1}{\em\ #2}{\ #3}}
\newcommand{\refer}[5]{{\sc #1}{\ #2}{\em\ #3}{\bf\ #4}{\ #5}}

\newtheorem{lem}{Lemma}[section]
\newtheorem{cor}[lem]{Corollary}
\newtheorem{teo}[lem]{Theorem}
\newtheorem{os}[lem]{Remark}
\newtheorem{defi}[lem]{Definition}
\newtheorem{prop}[lem]{Proposition}

\newcommand{\qed}{\thinspace\null\nobreak\hfill\hbox{\vbox{\kern-.2pt\hrule
 height.2pt depth.2pt\kern-.2pt\kern-.2pt \hbox to2.5mm{\kern-.2pt\vrule
 width.4pt \kern-.2pt\raise2.5mm\vbox to.2pt{}\lower0pt\vtop
 to.2pt{}\hfil\kern-.2pt \vrule
 width.4pt \kern-.2pt}\kern-.2pt\kern-.2pt\hrule height.2pt depth.2pt
 \kern-.2pt}}\par\medbreak}

\newcommand{\R}{\mathbb{R}}

\newcommand{\N}{\mathbb{N}}

\newcommand{\eps}{\varepsilon}

\newcommand{\ds}{\displaystyle}

\date{}
\begin{document}
\title{Heat Kernel estimates for some elliptic operators with unbounded diffusion coefficients}
\author{G. Metafune, C. Spina\thanks{Dipartimento di Matematica ``Ennio De
Giorgi'', Universit\`a del Salento, C.P.193, 73100, Lecce, Italy.
e-mail: giorgio.metafune@unisalento.it, chiara.spina@unisalento.it}}
\maketitle
\begin{abstract} We prove heat kernel bounds for the operator
$(1+|x|^\alpha)\Delta$ in $\R^N$, through Nash inequalities and
weighted Hardy inequalities.

\par


\end{abstract}

\section{Introduction and preliminary results}   \label{preliminari}
In this paper we prove heat kernel estimates for the operator
$$
L=m(x)(1+|x|^\alpha)\Delta
$$
in the whole space $\R^N$, under the assumption that $m$ is a
bounded, locally H\"older continuous function with $\inf m >0$. The
solvability of the elliptic and parabolic problem associated with
$L$, either in spaces of continuous functions or in $L^p$ spaces, has
been widely investigated in literature.  If $\alpha\leq 2$, $L$
generates an analytic semigroup both in $L^p(\R^N)$ and in
$C_0(\R^N)$, see \cite{for-lor}. If $\alpha>2$, the generation
results depend upon the space dimension $N$. If $N=1,\ 2$, $L$
generates a semigroup in $C_b(\R^N)$, the space of all continuous
and bounded functions on $\R^N$, but $C_0(\R^N)$ and $L^p(\R^N)$ are
not preserved. If $N\geq 3$, the resolvent and the semigroup map
$C_b(\R^N)$ into $C_0(\R^N)$ (see \cite{met-wack}) but $L^p(\R^N)$
is preserved if and only if $p>N/(N-2)$. For $N/(N-2)<p<\infty$ and
under the additional assumption that $m$ admits a finite limit at
infinity, the semigroup is also analytic in $L^p(\R^N)$ and, for
$m\equiv 1$, it is contractive if and only if $\alpha \le
(N-2)(p-1)$. We refer the reader to \cite{met-spina} for all these
results, as well as for domain characterization and spectral
properties of $L$. Here we only recall that the domain of $L$ in
$L^p(\R^N)$ coincides with the maximal one
$$D_{p,max}(L)=\{u \in W^{2,p}(\R^N): (1+|x|^\alpha)\Delta u \in L^p(\R^N)\}$$
and that the resolvents in $L^p(\R^N)$ and $L^q(\R^N)$ are
consistent, provided that $p,q >N/(N-2)$. Finally, the resolvent is
compact if and only if $\alpha >2$ and in this case the spectrum
consists of a sequence of negative eigenvalues $\lambda_n$ diverging
to $-\infty$.
\smallskip

Due to the local regularity of the coefficients, the semigroup
$(T(t))_{t\geq 0}$ generated by $L$ admits an integral kernel
$p(x,y,t)$, with respect to the Lebesgue measure (see e.g.
\cite{met-wack}), for which the following representation holds
$$T(t)f(x)=\int_{\R^N}p(x,y,t)f(y)dy.$$
However, the operator $L$ is symmetric with respect to the measure
$d\mu(x)=\left (1+|x|^\alpha\right)^{-1}dx$ and it is more
convenient to express $T(t)$ through a kernel with respect to
$d\mu$, namely
$$T(t)f(x)=\int_{\R^N}p_\mu(x,y,t)f(y)d\mu(y)$$
where
$$p_\mu(x,y,t)=(1+|y|^\alpha)p(x,y,t).$$
Our goal consists in obtaining upper bounds for the integral kernel
$p_\mu$ by working in  $L^2_\mu$ spaces and then deducing upper
bounds for $p$. This will be done by using the well-known
equivalence between Nash inequalities and ultracontractivity for
symmetric Markov semigroups, see \cite[Section 6.1]{ou}. \\

Throughout the paper the dimension $N$ will be always assumed to be
greater than or equal to $3$ and $\alpha$ will be a positive real
number.

We shall prove that for small $t$
$$p_\mu(x,y,t)\leq \frac{C}{t^\frac{N}{2}}(1+|x|^\alpha)^\frac{2-N}{4}(1+|y|^\alpha)^\frac{2-N}{4}$$
for $0<\alpha \le 4$, or
$$p_\mu(x,y,t)\leq \frac{C}{t^\frac{N+\alpha-2}{\alpha-2}}\phi(x)\phi(y)$$
for $2<\alpha\leq 4$ and
$$p_\mu(x,y,t)\leq \frac{C}{t^\frac{N}{2}}\phi(x)\phi(y)$$ for $\alpha\geq 4$.
Here $\phi$ is the first eigenfuncion of $L$ and satifies the bounds $C_1(1+|x|)^{2-N}\leq \phi(x)\leq C_2(1+|x|)^{2-N}$ for suitable $C_1, C_2>0$\\
We also show that the powers of $t$ appearing in the bounds above
are optimal. Estimates for large $t$ easily follow from the
semigroup law, since the semigroup dacays esponentially at infinity.
Observe that for $2<\alpha <4$ both the first and the second
estimate hold (see also Remark \ref{de}).

\subsection{Definition of $L$ via the quadratic form methods}
Consider the Hilbert spaces $L^2_\mu$, where $d\mu(x)=\left
(m(x)(1+|x|^\alpha)\right)^{-1}dx$, endowed with its canonical inner
product. Note that the measure $\mu$ is finite if and only if
 $\alpha >N$. Consider also the Sobolev space
 $$
 H=\{u \in L^2_\mu : \nabla u \in L^2\}
 $$
 endowed with the inner product
 $$
 (u,v)_H=\int_{\R^N}\left (u\bar{v}\, d\mu +\nabla u \cdot \nabla
 \bar{v}\, dx \right)
 $$
and let $\cal V$ be the closure of $C_c^1$ in $H$, with respect
to the norm of $H$. Observe that Sobolev inequality
\begin{equation} \label{sobolev}
\|u\|_{2^*}^2 \le C_2^2 \|\nabla u\|_2^2
\end{equation}
holds in $\cal V$ but not in $H$ (consider for example the case
where $\alpha >N$ and $u=1$). Here $2^*=2N/(N-2)$ and $C_2$ is the
best constant for which the equality above holds.\\
Next we introduce the continuous and weakly coercive symmetric form
\begin{equation} \label{forma}
a(u,v)= \int_{\R^N}\nabla u \cdot \nabla
 \bar{v}\, dx
\end{equation}
for $u,v \in \cal V$ and the self-adjoint operator $\cal L$
 defined by
 $$
 D({\cal L})=\{u\in L^2_\mu : {\rm there\ exists\ } f \in L^2_\mu:
 a(u,v)=-\int_{\R^N}f\bar{v}\, d\mu {\rm \ for\ every\  } v \in {\cal V} \} \qquad {\cal
 L}u=f.
$$
Since $a(u,u) \ge 0$, the operator ${\cal L}$ generates an analytic
semigroup of contractions $e^{t{\cal L}} $ in $L^2_\mu$. An
application of the Beurling-Deny criteria shows that the generated
semigroup is positive and $L^\infty$-contractive. For our purposes
we need  that the resolvents and the semigroups generated by ${\cal
L}$ and of $(L,D_{p,max}(L))$ are coherent. This is stated in the
following proposition. We refer to \cite[Proposition 7.4]{met-spina}
for its proof.

\begin{prop} \label{coerenza1}
$$
D({\cal L})\subset \{u \in {\cal V} \cap W^{2,2}_{loc}:
(1+|x|^\alpha) \Delta u \in L^2_\mu \}
$$
and ${\cal L} u=(1+|x|^\alpha)\Delta u$ for $u \in D({\cal L})$. If
$\lambda, t
>0$ and $f \in L^p \cap L^2_\mu$, then
$$
(\lambda-{\cal L})^{-1}f=(\lambda-L)^{-1}f
$$
and $$e^{t {\cal L}}f=T(t)f.
$$
\end{prop}
According with the above Proposition we write $T(t)$ for the
semigroup in $L^p$ with respect to the Lebesgue measure or in
$L^2_\mu$. It admits a positive integral kernel $p(x,y,t)$ with
respect to the Lebesgue measure, see \cite[Theorem 4.4]{met-wack},
for which the following representation holds
$$T(t)f(x)=\int_{\R^N}p(x,y,t)f(y)dy.$$
Clearly we have also
$$T(t)f(x)=\int_{\R^N}p_\mu(x,y,t)f(y)d\mu$$
with
$$p_\mu(x,y,t)=m(y)(1+|y|^\alpha)p(x,y,t).$$

\subsection{Eigenfunctions and eigenvalues of $L$ ($\alpha >2$)} Spectral
properties of $L$ have also been investigated in \cite[Section
7]{met-spina} where the following result has been proved. To unify
the notation, when $p=\infty$, $L^p$ stands for $C_0$. We recall
that the resolvent of $(L,D_{p,max}(L))$ is compact in $L^p$ if and
only if $\alpha>2$, see \cite{met-spina}, a condition that we assume
throughout this section.

\begin{prop} \label{indipendenza1}
If $N/(N-2)<p \le \infty$, $2<\alpha<\infty$, then the spectra of
$(L,D_{p,max}(L))$  and ${\cal L}$ coincide  and lie in $]-\infty,
0[$. They consist of a sequence $\lambda_n$ of eigenvalues, which
are simple poles of the resolvent and tend to $-\infty$. Each
eigenspace is finite dimensional and independent of $p$.
\end{prop}

The next propositions give a lower and upper bound of the first
eigenfunction of $L$. Better results will follow later from kernel
estimates.

\begin{prop}   \label{primautofunz}
Let $\lambda<0$ be the first eigenvalue of $L$ and $\phi$ be the
corresponding eigenfunction. Then there exists a positive constant
$C$ such that
$$\phi(x)\geq C(1+|x|)^{2-N}$$
for every $x\in\R^N$.
\end{prop}
{\sc Proof.} Since the kernel $p$ is positive, $T(t)$ is irreducible
and from \cite[Proposition 1.4.3]{davies1}) it follows that the
eigenspace relative to the first eigenalue is one-dimensional and
admists a strictly positive eigenfuncion $\phi$. Therefore
$c=\min_{B(1)} \phi>0$. Since $\phi \in D_{p,max}(L)$ for every
$N/(N-2)<p \le \infty$ and
$$\Delta \phi=\frac{\lambda\phi}{m(x)(1+|x|^\alpha)},$$
 from \cite[Section 4]{met-spina} we may
write
\begin{eqnarray*}
\phi(x)&=&\frac{1}{N(2-N)\omega_N}\int_{\R^N}\frac{\lambda\phi(y)}{m(y)|x-y|^{N-2}(1+|y|^\alpha)}dy
\\
&\geq &\frac{\lambda c}{N(2-N)\omega_N}
\int_{B(1)}\frac{1}{m(y)|x-y|^{N-2}(1+|y|^\alpha)}dy \\
&\geq & C\int_{B(1)}\frac{1}{|x-y|^{N-2}}dy \ge C (1+|x|)^{2-N}.
\end{eqnarray*}
\qed

Next we find upper bounds for the other eigenfunctions. Since the
spectrum is independent of $p$, the eigenfuncions $\phi_i$ belong to
$C_0$ and therefore are bounded.

\begin{prop}   \label{upper-autofunzioni}
Let $\phi_i$ be an eigenfunction of $L$  corresponding to the
eigenvalue $\lambda_i$. Then there exists a positive constant $C$
such that
$$|\phi_i(x)|\leq C|\lambda_i|^k\|\phi_i\|_\infty (1+|x|)^{2-N}$$
where $k=\ds\left[\frac{N}{\alpha-2}\right]+1$.
\end{prop}

The proof is based upon the following lemma whose proof can be found
in \cite[Lemma 6.1]{met-spina}.
\begin{lem} \label{comportamentoJ}
Let $\gamma,\ \beta>0$ such that $\gamma<N$ and $\gamma+\beta>N$.
Set
$$J(x)=\int_{\R^N}\frac{dy}{|x-y|^\gamma(1+|y|^\beta)}.$$
 Then $J$ is bounded in $\R^N$ and has the following  behaviour as $|x|$ goes to infinity
\begin{equation*}
J(x)\simeq\left\{
\begin{array}{ll}
c_1 |x|^{N-(\gamma+\beta)}& \textrm{if}\quad\beta<N \\
c_2|x|^{-\gamma}\log|x|   &\textrm{if}\quad\beta=N\\
c_3|x|^{-\gamma}  &\textrm{if}\quad\beta>N \end{array}\right.
\end{equation*}
for suitable positive constants $c_1, c_2, c_3$.
 \end{lem}

{\sc Proof} (Proposition \ref{upper-autofunzioni}). As in
Proposition \ref{primautofunz} we write
\begin{equation} \label{formula-autofunzione}
\phi_i(x)=\int_{\R^N}\frac{\lambda_i\phi_i(y)}{m(y)|x-y|^{N-2}(1+|y|^\alpha)}dy.
\end{equation}
It follows that
\begin{equation*}
|\phi_i(x)|\leq
C|\lambda_i|\|\phi_i\|_\infty\int_{\R^N}\frac{1}{|x-y|^{N-2}(1+|y|^\alpha)}dy.
\end{equation*}
By Lemma \ref{comportamentoJ},
\begin{equation*}
|\phi_i(x)|\leq\left\{
\begin{array}{ll}
c_1 |\lambda_i|\|\phi_i\|_\infty(1+|x|)^{2-\alpha}& \textrm{if}\quad\alpha<N \\
c_2|\lambda_i|\|\phi_i\|_\infty(1+|x|)^{2-N}\log|x|   &\textrm{if}\quad\alpha=N\\
c_3|\lambda_i|\|\phi_i\|_\infty(1+|x|)^{2-N}  &\textrm{if}\quad\alpha>N \end{array}\right.
\end{equation*}
If $\alpha>N$, we immediately deduce the claim. Otherwise, we
iterate the procedure as follows. Assuming $|\phi_i(x)|\leq c_1
|\lambda_i|\|\phi_i\|_\infty(1+|x|)^{2-\alpha}$ and inserting this
inequality in (\ref{formula-autofunzione}), we deduce
\begin{equation*}
|\phi_i(x)|\leq C|\lambda_i|^2\|\phi_i\|_\infty\int_{\R^N}\frac{(1+|y|)^{2-\alpha}}{|x-y|^{N-2}(1+|y|^\alpha)}dy
\leq C |\lambda_i|^2\|\phi_i\|_\infty\int_{\R^N}\frac{1}{|x-y|^{N-2}(1+|y|^{2\alpha-2})}dy
\end{equation*}
and, by Lemma \ref{comportamentoJ} again,  $|\phi_i(x)|\leq
C|\lambda_i|^2\|\phi_i\|_\infty(1+|x|)^{2(2-\alpha)}$. We iterate
this procedure $k$-times obtaining $|\phi_i(x)|\leq
C|\lambda_i|^k\|\phi_i\|_\infty(1+|x|)^{k(2-\alpha)}$ until
$k(\alpha-2) \le N$. The claim then follows since at the step
$k+1$.\qed

\section{Kernel estimates}  \label{ultracontr}
In this section we will  prove kernel estimates trough weighted Nash
inequalities involving suitable Lyapunov functions. The main tool
will be \cite[Theorem 2.5]{bakry} whose original formulation is due
to Wang, see \cite[Theorem 3.3]{wang}.

\begin{defi} \label{lyap}
 A Lyapunov function is a positive function $V$   such that
$$T(t)V(x)=\int_{\R^N}p_\mu(x,y,t) V(y)d \mu (y)\leq e^{ct}V(x)$$
for every $x\in\R^N$, $t>0$ and some positive constant $c$ (called a
Lyapunov constant).
\end{defi}

By following \cite[Definition 2.2]{bakry} we introduce weighted Nash
inequalities.

\begin{defi}  \label{weightedNash}
Let $V$ be a positive function on $\R^N$  and $\psi$ be a positive
function defined on $(0,\infty)$ such that $\ds\frac{\psi(x)}{x}$ is
increasing. A Dirichlet form $a$ on $L^2_\mu$ satisfies a weighted
Nash inequality with weight $V$ and rate function $\psi$ if
$$\psi\left(\frac{\|u\|^2_{L^2_\mu}}{\|uV\|^2_{L^1_\mu}}\right)\leq\frac{a(u,u)}{\|uV\|^2_{L^1_\mu}}$$
for all functions $u$ in the domain of the Dirichlet form such that $\|u\|^2_{L^2_\mu}>0$ and $\|uV\|^2_{L^1_\mu}<\infty$.
\end{defi}

 Next we state \cite[Theorem 2.5]{bakry}, adapted to our situation.
\begin{teo}  \label{bakry}
Let $(T(t))_{t\geq 0}$ be a symmetric Markov semigroup on $L^2_\mu$
with generator $L$. Assume that there exists a Lyapunov function $V$
with Lyapunov constant $c\geq 0$ and that the Dirichlet form
associated to $L$ satisfies a weighted Nash inequality with weight
$V$  and rate function $\psi$, integrable near infinity and not
integrable near zero.  Then
$$\|T(t)f\|_{L^2_\mu}\leq K(2t)e^{ct}\|fV\|_{L^1_\mu}$$ for all functions $f\in L^2_\mu$ such that $fV \in L^1_\mu$.
Here the function $K$ is defined by
\begin{displaymath}
K(t)=\sqrt{U^{-1}(t)}
\end{displaymath}
 where
$$U(t)=\int_t^\infty\frac{1}{\psi(u)}du.$$
\end{teo}
From \cite[Proposition 2.1]{bakry}, the following corollary follows.
\begin{cor}
Under the assumptions ons of Theorem \ref{bakry}, $(T(t))_{t\geq 0}$
has a kernel $p_\mu$ which satisfies
$$p_\mu(x,y,t)\leq K(t)^2e^{ct}V(x)V(y)$$ for all $t>0$, $x,\ y\in\R^N\times\R^N$.
\end{cor}

\begin{os}
In \cite[Theorem 2.5]{bakry},  $V$ is required to be in the domain
of the operator and to satisfy $LV\leq cV$. By the proof, however,  it is evident that such an
assumption is only used to ensure the validity of the inequality
\begin{equation} \label{uso-Lyapunov}
\int_{\R^N}VT(t)fd\mu\leq e^{ct}\int_{\R^N}Vfd\mu
\end{equation}
for all positive functions $f\in L^2_\mu$ such that $fV \in
L^1_\mu$.  Observe that,
if $V$ is a Lyapunov function according to our Definition
(\ref{lyap}), then, for every positive $f$, then by Fubini theorem
and the symmetry of $p_\mu$
\begin{eqnarray*}
\int_{\R^N}VT(t)fd\mu&=&\int_{\R^N}\int_{\R^N}p_\mu(t,x,y)V(x)f(y)d\mu(x)d\mu(y)\\
&=&\int_{\R^N}\int_{\R^N}p_\mu(t,x,y)V(y)f(x)d\mu(x)d\mu(y)=
\int_{\R^N}T(t)V\,fd\mu\leq e^{ct}\int_{\R^N}V\,fd\mu
\end{eqnarray*}
and hence (\ref{uso-Lyapunov}) holds.
\end{os}

\subsection{Intrinsic ultracontractivity}
We show kernel estimates through weighted Nash inequalities with
respect to the the first eigenfunction of $L$. More precisely we
will prove that Definition (\ref{weightedNash}) holds with Lyapunov
function $V$ given by the first eigenfunction $\phi$ of $L$ and rate
functions
$$\psi(t)=t^{1+\frac{2}{N}}\quad\textrm{or}\quad\psi(t)=t^{1+\frac{\alpha-2}{N+\alpha-4}}.$$

\begin{teo}  \label{nucleo-tpiccoli}
Assume that $\alpha >2$. Then the  kernel $p_\mu$ of the semigroup
generated by $ L$ satisfies
\begin{align*}
& p_\mu(x.y,t)\leq \frac{C}{t^\frac{N+\alpha-4}{\alpha-2}}(1+|x|)^{2-N}(1+|y|)^{2-N}\quad\quad 2< \alpha\leq 4\\
& p_\mu(x.y,t)\leq \frac{C}{t^\frac{N}{2}}(1+|x|)^{2-N}(1+|y|)^{2-N} \quad\quad \alpha\geq 4
\end{align*}
or, equivalently,
\begin{align*}
& p_\mu(x.y,t)\leq \frac{C}{t^\frac{N+\alpha-4}{\alpha-2}}\phi(x)\phi(y)\quad\quad 2< \alpha\leq 4\\
& p_\mu(x.y,t)\leq \frac{C}{t^\frac{N}{2}}\phi(x)\phi(y) \quad\quad
\alpha\geq 4
\end{align*}
 for every $0<t\leq 1$, $x,\ y\in\R^N$.
\end{teo}
{\sc Proof.} Let $u\in{\cal V}$ such that
$\|u\phi\|_{L^1_\mu}<\infty$. We treat separately the cases
$2<\alpha\leq 4$ and $\alpha\geq 4$. \\
{\bf Case $2<\alpha\leq 4$.} Set
$\theta=\frac{2(N+\alpha-4)}{\alpha-2}>2$. By H\"{o}lder's
inequality (with $p=(\theta+2)/(\theta-2)$, we obtain
\begin{eqnarray*}
\int_{\R^N}|u|^2 d\mu&=& \int_{\R^N}|u|^{\frac{2\theta}{\theta+2}+
\frac{4}{\theta+2}}
\frac{\phi^\frac{4}{\theta+2}}{\phi^\frac{4}{\theta+2}}d\mu
\\&\leq&\left(\int_{\R^N}|u|^\frac{2\theta}{\theta-2} \phi^\frac{4}{2-\theta}d\mu\right)^\frac{\theta-2}{\theta+2}
\left(\int_{\R^N}|u| \phi d\mu\right)^\frac{4}{\theta+2}
\end{eqnarray*}
and,  by recalling that
\begin{equation}\label{primautofunzione}
C_1(1+|x|)^{2-N}\leq \phi(x)\leq C_2(1+|x|)^{2-N},
\end{equation}
\begin{equation*}
\left(\int_{\R^N}|u|^2d\mu\right)^{1+\frac{2}{\theta}}\leq
C\left(\int_{\R^N}|u|^\frac{2\theta}{\theta-2}
(1+|x|)^{\frac{4(N-2)}{\theta-2}-\alpha} dx \right)^\frac{\theta-2}{\theta}
\left(\int_{\R^N}|u| \phi d\mu\right)^\frac{4}{\theta}.
\end{equation*}
Set now
\begin{align*}
& q=\frac{2\theta}{\theta-2}=\frac{2(N+\alpha-4)}{N-2},\qquad \beta=\frac{(\alpha-4)(N-2)}{2(N+\alpha-4)},\\
& \gamma=0,\qquad p=2,\qquad \nu=-\alpha
\end{align*}
and observe that, for $2<\alpha\leq 4$, the assumptions in Proposition \ref{weight-sobolev} are satisfied. Therefore
\begin{align*}
\left(\int_{\R^N}|u|^2 d\mu\right)^{1+\frac{2}{\theta}}\leq & C\left(\int_{\R^N}|\nabla u|^2 dx+\int_{\R^N}|u|^2 d\mu\right)\times\\
&\left(\int_{\R^N}|u| \phi d\mu\right)^\frac{4}{\theta}
\end{align*}
or, equivalently,
$$\|u\|_{L^2_\mu}^{2+\frac{4}{\theta}}\leq C \tilde{a}(u,u)\|u\phi \|_{L^1_\mu}^{\frac{4}{\theta}}$$
where
$$\tilde{a}(u,u)=\int_{\R^N}|\nabla u|^2 dx+\int_{\R^N}|u|^2 d\mu$$
is the quadratic form associated with the operator $L+I$.
By Theorem \ref{bakry}, we obtain
$$\tilde{p}_\mu(x,y,t)\leq \frac{C}{t^\frac{N+\alpha-4}{\alpha-2}}\phi(x)\phi(y)$$
for every $x,\ y\in\R^N$, $0<t\leq 1$ and where
$\tilde{p}_\mu=e^{-t}p_\mu$ is the kernel associated with  $L+I$. It
follows
$$p_\mu(x,y,t)\leq \frac{C}{t^\frac{N+\alpha-4}{\alpha-2}}\phi(x)\phi(y)$$
for every $x,\ y\in\R^N$, $0<t\leq 1$.\\
{\bf Case $\alpha\geq 4$.} As before, by H\"{o}lder's inequality
with $p=(N+2)/(N-2)$,
\begin{align*}
\left(\int_{\R^N}|u|^2d\mu\right)^{1+\frac{2}{N}}&\leq
\left(\int_{\R^N}|u|^\frac{2N}{N-2}
\phi^\frac{4}{2-N}d\mu\right)^\frac{N-2}{N}\left(\int_{\R^N}|u| \phi
d\mu\right)^\frac{4}{N}\\&\leq C \left(\int_{\R^N}|u|^\frac{2N}{N-2}
(1+|x|)^{4}\frac{1}{1+|x|^\alpha}dx\right)^\frac{N-2}{N}
\left(\int_{\R^N}|u| (1+|x|)^{2-N}d\mu\right)^\frac{4}{N}\\&\leq
C\left(\int_{\R^N}|u|^\frac{2N}{N-2}
dx\right)^\frac{N-2}{N}\left(\int_{\R^N}|u|
(1+|x|)^{2-N}d\mu\right)^\frac{4}{N}.
\end{align*}
 By Sobolev embedding \begin{align*}
\left(\int_{\R^N}|u|^2 d\mu\right)^{1+\frac{2}{N}}\leq &
C\left(\int_{\R^N}|\nabla u|^2 dx\right)\left(\int_{\R^N}|u| \phi
d\mu\right)^\frac{4}{N}
\end{align*}
and the kernel estimates follow as in the first part of the proof.
\qed

Large time estimates follow from the semigroup law and small times
estimates.
\begin{prop}  \label{tempi-grandi}
There exists a positive constant $C$ such that
\begin{equation*}
p_\mu(x,y,t)\leq Ce^{\lambda t}\phi(x)\phi(y)
\end{equation*}
for every $x,\ y\in\R^N$, $t\geq 1$, where $\phi$ is the first
eigenfunction of $L$ and $\lambda <0$ is the first eigenvalue.
\end{prop}
{\sc Proof.} Let $t\geq 1$. By the semigroup law,
$$p_\mu(x,y,t)=\int_{\R^N}p_\mu(x,z,t-1)p_\mu(z,y,1)d\mu(z).$$
By Theorem \ref{nucleo-tpiccoli}, we obtain
\begin{eqnarray*}
p_\mu(x,y,t)&\leq& C\int_{\R^N}p_\mu(x,z,t-1)\phi(z)\phi(y)d\mu(z)=C\phi(y)\int_{\R^N}p_\mu(x,z,t-1)\phi(z)d\mu (z)\\
&=&C\phi(y)e^{(t-1) L}\phi(x)=C_1 e^{\lambda t}\phi(x)\phi(y).
\end{eqnarray*}
\qed

\begin{os}
By recalling that  the kernel $p$ with respect to the Lebesgue
measure $dx$ satisfies $p_\mu(x,y,t)=(1+|y|^\alpha)p(x,y,t)$, we can
reformulate Theorem \ref{nucleo-tpiccoli} as follows
\begin{align*}
& p(x,y,t)\leq \frac{C}{t^\frac{N}{2}}(1+|x|)^{2-N}(1+|y|)^{2-N-\alpha} \quad\quad \alpha\geq 4 \\
& p(x,y,t)\leq
\frac{C}{t^\frac{N+\alpha-4}{\alpha-2}}(1+|x|)^{2-N}(1+|y|)^{2-N-\alpha}\quad\quad
2< \alpha\leq 4
\end{align*}
for $0 <t \le 1$. A similar remark holds also for Proposition
\ref{tempi-grandi}.
\end{os}

\begin{os}{\rm
The estimate
$$p_\mu(x,y,t)\leq
\frac{C}{t^\frac{N}{2}}(1+|x|)^{2-N}(1+|y|)^{2-N},
$$ which holds for $\alpha \ge 4$, is optimal among the estimates
of the form $c(t)\psi(x)\psi(y)$, since the space profile is that of
the first eigenfunction and the factor $t^{-N/2}$ cannot be improved
(by the local arguments of the proof of Theorem \ref{K}). Concerning
the estimate
$$
p_\mu(x,y,t)\leq
\frac{C}{t^\frac{N+\alpha-4}{\alpha-2}}(1+|x|)^{2-N}(1+|y|)^{2-N}
$$
which holds for $2<\alpha \le 4$, the same argument as before shows
its optimality with respect to the space variable. The optimality
with respect to time (among all estimates of the form
$t^{-\beta}(1+|x|)^{2-N}(1+|y|)^{2-N} $) is proved by the following
argument. Assume that
$$
p_\mu(x,y,t)\leq \frac{C}{t^\beta}(1+|x|)^{2-N}(1+|y|)^{2-N}
$$
holds for $\beta>0$ and $0<t \le 1$. By the argument of Proposition
\ref{tempi-grandi} it holds for every $t>0$ and then by
\cite[Theorem 2.10]{bakry} the weighted Nash inequality
\begin{equation}  \label{nash2}
\left(\int_{\R^N}|u|^2d\mu\right)^{1+\frac{1}{\beta}}\leq  C
a(u,u)\left(\int_{\R^N}|u|(1+|x|)^{2-N}d\mu\right)^\frac{2}{\beta}
\end{equation}

is valid  for every $u\in\cal V$. By replacing $u(x)$ with
$u(\lambda x)$, $\lambda>0$, in (\ref{nash2}) we obtain
$$\lambda^{\alpha-2+\frac{4-N-\alpha}{\beta}}\left(\int_{\R^N}\frac{|u|^2}{\lambda^\alpha+|x|^\alpha}dx\right)^{1+\frac{1}{\beta}}
\leq C
a(u,u)\left(\int_{\R^N}\frac{|u|}{\lambda^{\alpha+N-2}+|x|^{\alpha+N-2}}dx\right)^\frac{2}{\beta}$$
for every $u\in\cal V$ and $\lambda>0$. Letting $\lambda$ to zero,
such an inequality can be true only if the exponent of $\lambda$ is
nonnegative, that is if $\beta \ge (N+\alpha-4)/(\alpha-2)$. Observe
also that $(N+\alpha-4)/(\alpha-2) > N/2$ if and only if
$2<\alpha<4$.}
\end{os}

\subsection{Kernel estimates for $\alpha \le 4$}
In this subsection we find different kernel estimates which
hold for $0<\alpha \le 4$. In the case where $2<\alpha <4$, these
estimates are better than those of the preceeding subsection with
respect to the time variable, but worse with respect to the space
variables. We emphasize that the case $0<\alpha\leq 2$ was not
covered by the previous computations.

In the following proposition we show weighted Nash inequalities with
respect to the weight $V=(1+|x|^\alpha)^\frac{2-N}{4}$.

\begin{prop} \label{nash}
Let $u\in\cal V$. Then
$$\left(\int_{\R^N}|u|^2d\mu\right)^{1+\frac{2}{N}}\leq a(u,u)\left(\int_{\R^N}|u|(1+|x|^\alpha)^{\frac{2-N}{4}}d\mu\right)^\frac{4}{N}.$$
\end{prop}
{\sc Proof.} Let $u\in\cal V$. Then, by H\"{o}lder's inequality,
\begin{eqnarray*}
\int_{\R^N}|u|^2d\mu & \le &C\int_{\R^N}|u|^2(1+|x|^\alpha)^{-1}dx=C\int_{\R^N}|u|^{\frac{2N}{N+2}}|u|^\frac{4}{N+2}(1+|x|^\alpha)^{-1}dx\\
&\leq&
C\left(\int_{\R^N}|u|^{\frac{2N}{N+2}\frac{N+2}{N-2}}dx\right)^\frac{N-2}{N+2}\left(|u|(1+|x|^\alpha)^{-\frac{N+2}{4}}dx\right)^\frac{4}{N+2}\\
&=&C\left(\int_{\R^N}|u|^{2^*}dx\right)^{\frac{1}{2^*}\frac{2N}{N+2}}\left(\int_{\R^N}|u|(1+|x|^\alpha)^{\frac{2-N}{4}}d\mu\right)^\frac{4}{N+2},
\end{eqnarray*}
where $2^*=2N/(N-2)$.  By Sobolev embedding,
$$\int_{\R^N}|u|^2d\mu\leq C\left(\int_{\R^N}|\nabla u|^{2}dx\right)^\frac{N}{N+2}
\left(\int_{\R^N}|u|(1+|x|^\alpha)^{\frac{2-N}{4}}d\mu\right)^\frac{4}{N+2}$$
hence
$$\left(\int_{\R^N}|u|^2d\mu\right)^{1+\frac{2}{N}}\leq  C a(u,u)\left(\int_{\R^N}|u|(1+|x|^\alpha)^{\frac{2-N}{4}}d\mu\right)^\frac{4}{N}.$$
\qed

 Next we show that $V$ is a Lyapunov function.
\begin{lem}  \label{lyapunov}
The function $V(x)=(1+|x|^\alpha)^{\beta}$ satisfies $LV(x)=m(x)(1+|x|^\alpha)\Delta V(x)\leq 0$ for $\frac{2-N}{\alpha}\leq \beta<0$.
\end{lem}
{\sc Proof.} In fact
\begin{eqnarray*}
(1+|x|^\alpha)\Delta V&=&\alpha\beta
|x|^{\alpha-2}V(x)\left[\alpha(\beta-1)\frac{|x|^\alpha}{1+|x|^\alpha}+\alpha-2+N\right]\\&\leq&
\alpha\beta(\alpha\beta-2+N)|x|^{\alpha-2}V \le 0
\end{eqnarray*}
if $\frac{2-N}{\alpha}\leq \beta<0$.
\qed
\begin{cor} \label{weighted1}
If $\alpha\leq 4$, the function $V(x)=(1+|x|^\alpha)^{\frac{2-N}{4}}$ satisfies $LV\leq 0$.
\end{cor}

\begin{lem} \label{non-Lyapunov}
 If $0 < \alpha \le 4$, then $V(x)=(1+|x|^\alpha)^\frac{2-N}{4}$ is a Lyapunov functions for $L$ with Lyapunov constant $c=0$.
\end{lem}
{\sc Proof.} Observe that $V\in C_0(\R^N)$. Let $\lambda>0$ and set
$u=R(\lambda, L)V$, where the resolvent is understood  in
$C_0(\R^N)$, see \cite{met-wack} or \cite[Section 6]{met-spina}.
Since $LV \le 0$ then
$$\lambda\left(\frac{V}{\lambda}-u\right)-L\left(\frac{V}{\lambda}-u\right)\geq
0.$$By the maximum principle it follows that $V\geq \lambda
R(\lambda,L)V$. By iteration the last inequality implies
$$V\geq \lambda^n R(\lambda,L)^n V$$
for every $n\in\N$. Then
$$T(t)V=\lim_{n\to\infty}\left[\frac{n}{t}R\left(\frac{n}{t}, L\right)V\right]^n\leq V$$ in
$C_0(\R^N)$. \qed

From Proposition \ref{nash}, Lemma \ref{non-Lyapunov} and
\cite[Corollary 2.8]{bakry} the following result follows.

\begin{teo}  \label{stime-fin}
If $0 <\alpha \le 4$ the semigroup generated by $ L$  has a kernel
$p_\mu$ with respect to the measure $d\mu$ that satisfies the
following bounds
$$
 p_\mu(x,y,t)\leq
\frac{C}{t^\frac{N}{2}}(1+|x|^\alpha)^{\frac{2-N}{4}}(1+|y|^\alpha)^{\frac{2-N}{4}}
$$ for every $t>0$, $x,\ y\in\R^N$.
\end{teo}

\begin{os}
By recalling that  the kernel $p$ with respect to the Lebesgue
measure $dx$ satisfies $p_\mu(x,y,t)=(1+|y|^\alpha)p(x,y,t)$, we
deduce for $0<\alpha \le 4$
$$
 p(x,y,t)\leq \frac{C}{t^\frac{N}{2}}(1+|x|^\alpha)^{\frac{2-N}{4}}(1+|y|^\alpha)^{\frac{2-N}{4}-1}
 $$
 for every $t>0$, $x,\ y\in\R^N$.
\end{os}

\begin{os}  \label {de}
When $2 < \alpha<4$, writing $p_\mu=p_\mu^\theta p_\mu^{1-\theta}$ for every $0 \le \theta \le 1$, we can combine the estimates of Theorems \ref{nucleo-tpiccoli}, \ref{stime-fin} thus obtaining a family of bounds depending on a parameter $\theta$.
\end{os}

\section{Some consequences}
In this section we assume $\alpha >2$ and deduce further properties
of the eigenvalues and eigenfunctions of $L$ from kernel estimates.
Let us denote by
$$\lambda_1\geq\lambda_2\geq...$$ the eigenvalues of $L$, repeated according to their multiplicity, and by $\phi_n$ the corresponding eigenfunctions,
which we assume to be normalized in $L^2_\mu$. We write $\phi$ for
$\phi_1$.

\begin{prop}
The following estimates hold.
$$|\phi_n(x)|\leq C|\lambda_n|^{\frac{N}{4}}|\phi(x)|$$ for every $x\in \R^N$ and for $\alpha\geq
4$.
$$|\phi_n(x)|\leq C|\lambda_n|^{\frac{N+\alpha-4}{2(\alpha-2)}}|\phi(x)|$$ and
$$|\phi_n(x)|\leq C|\lambda_n|^{\frac{N}{4}}(1+|x|^\alpha)^\frac{2-N}{4}$$ for every $x\in \R^N$ and $2<\alpha\leq 4$.
\end{prop}
{\sc Proof.}  Let $\tilde{T}(t)f=\frac{1}{V}T(t)(Vf)$ with $V=\phi$
or $V=(1+|x|^\alpha)^\frac{2-N}{4}$. Assume e.g. $\alpha \ge 4$.
Theorem \ref{nucleo-tpiccoli} and Proposition \ref{tempi-grandi}
show that $\tilde{T}(t)$ is a bounded by $Ct^{-N/2}$ from $(L^1,
V^2d\mu)$ into $L^\infty$ and (\ref{uso-Lyapunov}) shows that it is
a contraction from $(L^1, V^2d\mu)$ into itself. Then $\tilde{T}(t)$
is  bounded by $Ct^{-N/4}$ from $(L^2, V^2d\mu)$ into $L^\infty$
and hence
\begin{eqnarray*}
\left\|\tilde{T}(t)\left(\frac{\phi_n}{\phi}\right)\right\|_\infty\leq
\frac{C}{t^\frac{N}{4}}\left\|\frac{\phi_n}{\phi}\right\|_{(L^2,\phi^2d\mu)}=\frac{C}{t^\frac{N}{4}}.
\end{eqnarray*}
Since
$$\tilde{T}(t)\left(\frac{\phi_n}{\phi}\right)=\frac{1}{\phi}T(t)\phi_n=\frac{e^{\lambda_nt}\phi_n}{\phi},$$
it follows that
$$\left|\frac{\phi_n(x)}{\phi(x)}\right|\leq \frac{C e^{-\lambda_n t}}{t^\frac{N}{4}}$$
for every $x\in\R^N$ and $t>0$. Minimizing over $t$ we obtain
$$\left|\phi_n(x)\right|\leq C |\lambda_n|^\frac{N}{4}|\phi(x)|$$ for every $x\in\R^N$.
The estimates for $2<\alpha\leq 4$ follow in a similar way. \qed

For $\lambda>0$ let $N(\lambda)$ be the number of $\lambda_j$ such
that $|\lambda_j|\leq\lambda$. The kernel estimates allow to deduce
some information on the distribution of the eigenvalues. The
following result is usually obtained as a corollary of the classical
Mercer's Theorem. We refer to \cite[Proposition 4.1]{met-spina2} for
a simple proof based on the semigroup property. The convergence of
the integral below is easily verified using Theorem
\ref{nucleo-tpiccoli} for $\alpha \ge 4$ and Theorem \ref{stime-fin}
for $2<\alpha <4$.
\begin{prop}   \label{stimadiag}
Let $t>0$. Then
$$\sum_{n=1}^\infty e^{\lambda_nt}=\int_{\R^N}p_\mu(x,x,t)\, d\mu(x) <\infty.$$
\end{prop}
 The following proposition, which is a weaker (and
elementary) version of Karamata's Theorem, allows to deduce
informations on $N(\lambda)$. For its proof we refer to
\cite[Proposition 7.2]{met-spina2}.

\begin{prop}  \label{Karamata-weak}
Let $r>0,\ C_1>0$ such that
\begin{equation}   \label{limsup}
\limsup_{t\to 0}t^r\sum_{n\in \N}e^{\lambda_n t}\leq C_1.
\end{equation}
Then
$$\limsup_{\lambda\to\infty}\lambda^{-r}N(\lambda)\leq C_1\frac{e^r}{r^r}.$$
Moreover if (\ref{limsup}) holds and
\begin{equation}   \label{liminf}
\liminf_{t\to 0}t^r\sum_{n\in\N}e^{\lambda_n t}\geq C_2
\end{equation}
for some $C_2>0$ then
$$\liminf_{\lambda\to\infty}\lambda^{-r}N(\lambda)\geq C_3$$ for some positive $C_3$.
\end{prop}

\begin{teo} \label{K}
Let $N(\lambda)$ be defined as before. Then
$$\limsup_{\lambda\to\infty}\lambda^{-\frac{N}{2}}N(\lambda)\leq C_1$$ and
$$\liminf_{\lambda\to\infty}\lambda^{-\frac{N}{2}}N(\lambda)\geq C_2$$ for some positive $C_1,\ C_2$.
\end{teo}
{\sc Proof.} By Theorems \ref{nucleo-tpiccoli}, \ref{stime-fin} we
deduce
$$\limsup_{\lambda\to\infty}\lambda^{-\frac{N}{2}}N(\lambda)\leq C_1$$ for some positive constant $C_1$. On the other
hand, for $|x| \le 1$  let $\Omega=B(x,1)$ and denote by $p_\Omega$
the heat kernel of the restriction of $L$ to $\Omega$ with Dirichlet
boundary conditions. Since we have (see  \cite[Remark 6]{met-ouh}
$$p_\mu(x,x,t)\geq p_\Omega(x,x,t)\geq \frac{C}{t^{\frac{N}{2}}}$$ for $0<t\leq 1$ and constant $C$  independent of $|x| \le 1$, it follows that
$$ \int_{\R^N}p_\mu(x,x,t)d\mu(x)\ge \int_{B(0,1)}p_\mu(x,x,t)d\mu(x)\ge \frac{C}{t^{\frac{N}{2}}}$$
for $0<t \le 1$. Proposition \ref{Karamata-weak} implies that
$$\liminf_{\lambda\to\infty}\lambda^{-\frac{N}{2}}N(\lambda)\geq C_2$$ for some positive $C_2$.\qed

\section*{Appendix: some weighted Sobolev inequalities}

The weighted Sobolev inequalities below have been a major tool to
prove kernel estimates. For reader's convenience  we provide a
direct proof following the methods of \cite{bazan}). We point out
that such type of inequalities  follow also by more general results
about Sobolev inequalities with respect to weights satisfying
Muckenhoupt-type conditions, sse \cite{muck}, \cite{muck1}.

\begin{prop}  \label{weight-sobolev}
Let  $\beta,\ \gamma,\ \nu\in \R,\ \gamma-1\leq \beta\leq\gamma$,
$1< p \le q<\infty$  satisfy
$$0\leq\frac{1}{p}-\frac{1}{q}=\frac{1-\gamma+\beta}{N}, \quad\quad
N+p(\gamma-1)\neq 0, \quad \quad p\leq q\leq p^*, \quad\quad p<N.$$
Then there exists a positive constant $C$ such that for every $u\in
C_c^\infty(\R^N)$
\begin{equation}  \label{weigted-Sobolev}
\left(\int_{\R^N}(1+|x|)^{q\beta}|u(x)|^qdx\right)^\frac{1}{q}\leq C\left(\int_{\R^N}(1+|x|)^{\gamma p}|\nabla u(x)|^pdx\right)^\frac{1}{p}+ C\left(\int_{\R^N}(1+|x|)^{\nu}|u(x)|^pdx\right)^\frac{1}{p}.
\end{equation}

\end{prop}

First we prove an homogeneous version of the previous inequality.
\begin{lem}
Let  $\beta,\ \gamma\in \R,\ \gamma-1\leq \beta\leq\gamma$,  $1< p
\le q<\infty$  satisfy
$$0\leq\frac{1}{p}-\frac{1}{q}=\frac{1-\gamma+\beta}{N},\quad\quad
N+p(\gamma-1)\neq 0,   \quad \quad p\leq q\leq p^*, \quad\quad
p<N.$$ Then there exists a positive constant $C$ such that for every
$u\in C_c^\infty(\R^N\setminus\{0\})$
 \begin{equation}  \label{fuori-palla}
\left(\int_{\R^N}|x|^{q\beta}|u(x)|^qdx\right)^\frac{1}{q}\leq C\left(\int_{\R^N}|x|^{p\gamma}|\nabla u(x)|^pdx\right)^\frac{1}{p}.
\end{equation}
\end{lem}

{\sc Proof.} We follow \cite{bazan}. By applying the divergence
theorem to suitable vector fields, we prove (\ref{fuori-palla}) in
correspondence of $\beta=\gamma-1$ and $p=q$, then by classical
embeddings theorems we prove (\ref{fuori-palla}) for $\gamma=\beta$
and $q=p^*$ and finally, by interpolation, we deduce
(\ref{fuori-palla}) in the general case. Let $u\in
C_c^\infty(\R^N\setminus\{0\})$. Since
$${\rm div} \left(\frac{x|x|^{(\gamma-1)p}}{N+p(\gamma-1)}\right)=|x|^{(\gamma-1)p}$$ by the divergence theorem and H\"{o}lder's inequality,
\begin{eqnarray*}
& &\int_{\R^N}|x|^{(\gamma-1)p}|u(x)|^p dx =\int_{\R^N}{\rm div} \left(\frac{x|x|^{(\gamma-1)p}}{N+p(\gamma-1)}\right)|u(x)|^p dx=\\
&-&\int_{\R^N}\frac{x|x|^{(\gamma-1)p}}{N+p(\gamma-1)}\nabla(|u(x)|^p)\leq
\frac{p}{N+p(\gamma-1)}\int_{\R^N}|x|^{(\gamma-1)(p-1)+\gamma}|u(x)|^{p-1}|\nabla
u(x)|dx\\& \leq&
C\left(\int_{\R^N}|x|^{(\gamma-1)p}|u(x)|^pdx\right)^\frac{1}{p^{'}}\left(\int_{\R^N}|x|^{\gamma
p}|\nabla u(x)|^pdx\right)^\frac{1}{p}.
\end{eqnarray*}
It follows
\begin{equation} \label{divergenza}
\left(\int_{\R^N}|x|^{(\gamma-1)p}|u(x)|^pdx\right)^\frac{1}{p}\leq C\left(\int_{\R^N}|x|^{\gamma p}|\nabla u(x)|^pdx\right)^\frac{1}{p}.
\end{equation}
Let us now consider the case $q=p^*$. Setting
$f(x)=|x|^{\gamma}u(x)$, by Sobolev inequality we have
\begin{eqnarray*}
& &\int_{\R^N}|x|^{\gamma p^*}|u(x)|^{p^*} dx=\int_{\R^N}|f(x)|^{p^*} dx\leq
 C\left(\int_{\R^N}|\nabla f(x)|^{p}dx\right)^\frac{p^*}{p}\\&\leq& C\left(\int_{\R^N}|x|^{(\gamma-1)p}|u(x)|^{p}dx\right)^\frac{p^*}{p}+C\left(\int_{\R^N}|x|^{\gamma p}|\nabla u(x)|^{p}dx\right)^\frac{p^*}{p}.
\end{eqnarray*}
Using (\ref{divergenza}) to estimate the first addendum in the right
hand side of the previous expression, we arrive at
\begin{equation}  \label{classicoSobolev}
\left(\int_{\R^N}|x|^{\gamma p^*}|u(x)|^{p^*} dx\right)^\frac{1}{p^*}\leq C\left(\int_{\R^N}|x|^{\gamma p}|\nabla u(x)|^pdx\right)^\frac{1}{p}.
\end{equation}
Let now $p\leq q \leq p^*$ and  $\gamma-1\leq \beta\leq \gamma$.
There exists $\theta\in [0,1]$ such that $\beta
=(1-\theta)(\gamma-1)+\theta\gamma$. By the assumption $\ds
\frac{1}{p}-\frac{1}{q}=\frac{1-\gamma+\beta}{N}$ it follows that
$q=\ds\frac{Np}{N-\theta p}=(1-c)p+cp^*$ with
$c=\ds\frac{\theta(N-p)}{N-\theta p}\in [0,1]$. Therefore,
H\"{o}lder's inequality yields
\begin{eqnarray*}
\int_{\R^N}|x|^{q\beta}|u(x)|^qdx&=&\int_{\R^N}|x|^{(1-c)(\gamma-1)p+\gamma c p^*}|u(x)|^{(1-c)p+cp^*}dx\\
&=& \int_{\R^N}\left(|x|^{(\gamma-1)p}|u(x)|^p\right)^{1-c}\left(|x|^{\gamma p^*}|u(x)|^{p^*}\right)^{c}dx\\&\leq&
\left(\int_{\R^N}|x|^{(\gamma-1)p}|u(x)|^p\right)^{1-c}
\left(\int_{\R^N}|x|^{\gamma p^*}|u(x)|^{p^*}\right)^c.
\end{eqnarray*}
By (\ref{divergenza}) and (\ref{classicoSobolev}), we obtain
$$\left(\int_{\R^N}|x|^{q\beta}|u(x)|^qdx\right)^\frac{1}{q}\leq C \left(\int_{\R^N}|x|^{\gamma p}|\nabla u(x)|^pdx\right)^\frac{1}{p}.$$
\qed

{\sc Proof.} (Proposition \ref{weight-sobolev}) Let $u\in
C_c^\infty(\R^N)$ and fix $\eta\in C_c^\infty(\R^N),\ 0\leq \eta\leq
1$ such that $\eta=1$ in $B(1)$ and $\eta=0$ in $\R^N\setminus
B(2)$. By Sobolev embedding,
\begin{eqnarray*}
\left(\int_{ B(2)}(1+|x|)^{q\beta}|\eta(x)u(x)|^qdx\right)^\frac{1}{q}&\leq& C\left(\int_{ B(2)}|\eta(x)u(x)|^qdx\right)^\frac{1}{q}\leq  C\left(\int_{ B(2)}|\nabla (\eta(x)u(x))|^p dx\right)^{\frac{1}{p}}\\&\leq& C\left(\int_{ B(2)}(1+|x|)^{\gamma p}|\nabla(\eta(x) u(x))|^p dx\right)^{\frac{1}{p}}.
\end{eqnarray*}
By combining with (\ref{fuori-palla}), we deduce
\begin{eqnarray*}
& &\left(\int_{\R^N}(1+|x|)^{q\beta}|u(x)|^qdx\right)^\frac{1}{q}=\left(\int_{\R^N}(1+|x|)^{q\beta}|\eta(x) u(x)+(1-\eta)(x)u(x)|^qdx\right)^\frac{1}{q}\\
&\leq& C\left(\int_{\R^N}(1+|x|)^{\gamma p}|\nabla
u(x)|^pdx\right)^\frac{1}{p}+ C\left(\int_{B(2)}(1+|x|)^{\gamma
p}|\nabla\eta(x)|^p|u(x)|^pdx\right)^\frac{1}{p}.
\end{eqnarray*}

\qed


\begin{thebibliography}{99}
%\bibitem{arendt}
%\refer{W. Arendt:} {Gaussian estimates and interpolation of the spectrum in $L^p$,}{Diff. Int. Eq.,}{Vol. 7,n. 5}{(1994), 1153-1168}.

\bibitem{bazan}
\refer{A. Bazan, W. Neves: }{The Kaffarelli-Kohn-Niremberg's Inequality for arbitrary norms}{preprint}{×}{×}.

\bibitem{bakry}
\refer{D. Bakry, F. Bolley, I. Gentil, P. Maheux:} {Weighted Nash Inequalities,} {arXiv: 1004.3456.} {} {}

\bibitem{davies1}
\refbook{E. B. Davies: }{One-Parameter Semigroups, }{ Academic
Press, 1980}

\bibitem{davies}
\refbook{E. B. Davies: }{Heat kernels and spectral theory, }{
Cambridge University Press, 1989}

\bibitem{engel-nagel}
\refbook{K.J. Engel, R. Nagel:} {One parameter semigroups for linear evolutions equations,} {Springer-Verlag, Berlin, (2000)}.

\bibitem{for-lor}
\refer{S. Fornaro, L. Lorenzi:} {Generation results for elliptic operators with unbounded diffusion coefficients in $L^p$ and $C_b$-spaces,} {Discrete and continuous dynamical sistems,} {Vol. 18, N.4} {(2007),  747-772}.

\bibitem{met-ouh}
\refer{G. Metafune, E.M. Ouhabaz, D. Pallara:} {Long time behavior of heat kernels of operators with unbounded drift terms, }{ to appear.}{} {}{}

\bibitem{met-spina}
\refer{G. Metafune, C. Spina:} {Elliptic operators with unbounded diffusion coefficients in $L^p$ spaces,} {to appear, Annali di Pisa}. {}

\bibitem{met-spina2}
\refer{G. Metafune, C. Spina:} {Kernel estimates for a class of  certain
Schr\"odinger semigroups, }{Journal of Evolution Equations, }{7} {(2007),}{ 719-742.}


\bibitem{met-wack}
\refer{G. Metafune, D. Pallara, M. Wacker:} {Feller Semigroups on
$\R^N$,} {Semigroup Forum,} {65} {(2002),  159-205}.

\bibitem{muck}
\refer{B. Muckenhoupt:}{Hardy's inequalities with weights,}{ Studia
Math.,} {34}{ (1972), 31-38}.

\bibitem{muck1}
{B. Muckenhoupt, R. Wheeden:}{ Weighted norm inequalities for
fractional integrals,}{Trans. Amer. Math. Soc.,}{192}{ (1974)
261-274}.

\bibitem{ou}
\refbook{E. M. Ouhabaz:} {Analysis of Heat Equations on Domains,}
{Princeton University Press}.


\bibitem{wang}
\refer{F. Y. Wang:} {Functional inequalities and spectrum estimates:
the infinite measure case,} {J. Funct. Anal.,} {194} {(2002),
288-310}.

\end{thebibliography}
\end{document}